\documentclass{svproc}
\usepackage{graphicx,graphics}
\usepackage{url}

\begin{document}
\mainmatter              
\title{Prabhakar discrete-time generalization of the time-fractional Poisson process and related random walks}
\titlerunning{Prabhakar discrete-time process}  
%
\author{Thomas M. Michelitsch\inst{1} \and Federico Polito\inst{2}
\and
Alejandro P. Riascos\inst{3}}
\authorrunning{T. M. Michelitsch, F. Polito \& A.P. Riascos} 
%
\tocauthor{}
\institute{Sorbonne Universit\'e, Institut Jean le Rond d’Alembert, 
CNRS UMR 7190,\\ 4 place Jussieu, 75252 Paris cedex 05, France\\
\email{michel@lmm.jussieu.fr}\\[1mm]
\and
Department of Mathematics ``Giuseppe Peano'', University of Torino, Italy \\[1mm]
\and
Instituto de F\'isica, Universidad Nacional Aut\'onoma de M\'exico,\\
Apartado Postal 20-364, 01000 Ciudad de M\'exico, M\'exico
}

\maketitle              

\begin{abstract}
In recent years a huge interdisciplinary field has emerged which is devoted to the `complex dynamics' of anomalous transport with long-time memory and non-markovian features.
It was found that the framework of fractional calculus and its generalizations are able to capture these phenomena.
Many of the classical models are based on continuous-time renewal processes and use the Montroll--Weiss continuous-time random walk (CTRW) approach. On the other hand their discrete-time counterparts are rarely considered in the literature despite their importance in various applications.
The goal of the present paper is to give a brief sketch of our recently introduced discrete-time Prabhakar generalization of the fractional Poisson process and the related discrete-time random walk (DTRW) model. We show that this counting process is connected with the continuous-time Prabhakar renewal process by a (`well-scaled') continuous-time limit.
We deduce the state probabilities and discrete-time generalized fractional Kolmogorov-Feller equations governing the Prabhakar DTRW and discuss
effects such as long-time memory (non-markovianity) as a hallmark of the `complexity' of the process.
\end{abstract}
\section{Introduction}
\label{Intro}
It is well-known that phenomena of anomalous transport and diffusion are not compatible with the classical picture of exponentially distributed waiting time patterns. Indeed it was found that random walks subordinated to an independent renewal process with fat-tailed waiting time densities (in other words Montroll--Weiss continuous-time random walk (CTRW) models formulated in the framework of fractional calculus) offer an appropriate description (see e.g.\ \cite{MetzlerKlafter2000,Meerschert-et-al2018,GorenfloMainardi2013,BeghinOrsinger2009,Laskin2003}). Further, these models have profound connections with semi-Markov processes \cite{Levy1956}. In order to obtain a greater flexibility to describe real-world situations, 
some generalized models have emerged \cite{Mi-PoRia-fractal-fract}, among them
the so-called Prabhakar generalization \cite{PolitoCahoy2013,GiustiPolitoMainardi-etal2020,TMM-APR-PhysicaA2020,MiRia2020}.
In most cases such as the Montroll--Weiss CTRW \cite{MontrollWeiss1965} the time variable is continuous. However, there are many systems in nature having intrinsically discrete time scales; such cases are rarely considered in the literature. Essential elements of this theory were developed recently \cite{Pachon-Polito-Ricciuti2020} (see also the references therein).

In the present letter we give a brief sketch of our recent model of a Prabhakar discrete-time generalization of the fractional Poisson process and its application to stochastic motions on undirected graphs. We refer the reader to our recent paper \cite{Michelitsch-Polito-Riascos2021} for details.
\section{The Prabahakar discrete-time counting process}
\label{PDTP}
We consider a counting process with integer arrival times characterized by the variable 
\begin{equation}
\label{PDTP_counting}
 {\cal J}_n = \sum_{j=1}^n Z_j, \hspace{0.5cm}  Z_j \in {\mathbf N}=\{1,2,\dots\}  , \hspace{0.5cm} {\cal J}_0=0,
\end{equation}
where $Z_j$ are IID non-zero integer inter-arrival times following the distribution $\theta_{\alpha}^{(\nu)}(t,\xi)$, $t \in {\mathbf N}$, $\nu>0$, $\alpha \in (0,1]$, $\xi>0$,  and having generating function
\begin{equation} 
\label{genfugen}
{\bar 
\theta}_{\alpha}^{(\nu)}(u)  = {\bar f}(u) {\bar \varphi}_{\alpha}^{(\nu)}(u) ,
\quad {\bar \varphi}_{\alpha}^{(\nu)}(u)= \frac{\xi^{\nu}}{(\xi+(1-u)^{\alpha})^{\nu}}.
\end{equation}
In this expression we have ${\bar f}(u) = u \sum_{k=1}^{\infty} f(k)u^{k-1}$ with ${\bar f}(u)|_{u=1}=1$ fulfilling
the desired initial condition
${\bar f}(u)|_{u=0} = f(t)\big|_{t=0}=0$, ensuring almost sure nonzero inter-arrival times. 
In (\ref{genfugen}) the parameter $\xi$ defines a characteristic time-scale and crucial here is that
${\bar f}(u)$ is free of the parameter $\xi$.

\begin{figure}[!t] 
\begin{center}
\includegraphics*[width=1.0\textwidth]{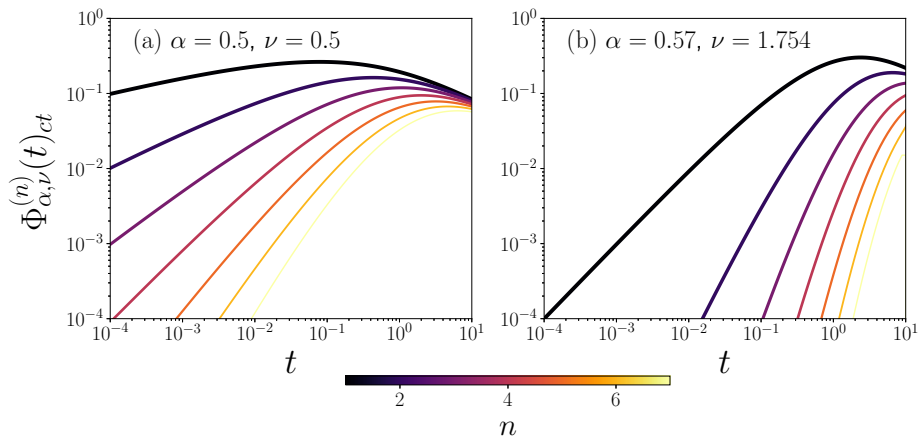}
\end{center}
\vspace{-5mm}
\caption{\label{Fig1} GFPP state probabilities $\Phi_{\alpha,\nu}^{(n)}(t)_{ct}$ of (\ref{well-sclaed_state-prob}) versus $t$ for different values of $n=1,2,\ldots,7$ and $\xi_0=1$. (a) $\alpha = 0.5$ and $\nu = 0.5$, (b) $\alpha = 0.57$ and $\nu = 1.754$.  }
\end{figure}

Consider the simplest case: ${\bar f}(u)=u$ with ${\bar 
\theta}_{\alpha}^{(\nu)}(u)=u {\bar \varphi}_{\alpha}^{(\nu)}(u)$. In this case, we call the counting process `{\it Prabhakar discrete-time process (PDTP)}'. The PDTP stands out since it contains several classical cases: for $\nu =1$ and $0<\alpha<1$, it recovers the `{\it fractional Bernoulli process}'
with the `discrete-time Mittag--Leffler distribution' (of so-called `type A') DML$_A$  \cite{Pachon-Polito-Ricciuti2020}, and for $\nu=1$ and $\alpha=1$ the standard Bernoulli counting process.
The PDTP waiting time density then has the form $\theta_{\alpha}^{(\nu)}(t,\xi)=\varphi_{\alpha}^{(\nu)}(t-1,\xi)$ with
\begin{equation}
\label{discrete-Prabhakar}
\begin{array}{l}
\displaystyle
\varphi_{\alpha}^{(\nu )}(k,\xi)  \displaystyle =
\left\{\begin{array}{clc} \displaystyle 
\frac{\xi^{\nu}}{k!}\sum_{m=0}^{\infty} \frac{(-1)^m (\nu )_{m} \xi^m}{m!}  \frac{\Gamma[\alpha(m+\nu)+k]}{\Gamma[\alpha(m+\nu)]} , &\hspace{0.2cm} 0 < \xi < 1, \\ \\
\displaystyle \frac{(-1)^k}{k!}\sum_{m=0}^{\infty} \frac{(-1)^m (\nu )_{m} \xi^{-m}}{m!}
\frac{\Gamma[\alpha m+1]}{\Gamma[\alpha m-k+1]} , & \hspace{0.3cm}\xi > 1,
\end{array}\right. 
\end{array} 
\end{equation}
$k \in {\mathbf N_0} = \{ 0,1,2,\ldots \}$ and where we employed the Pochhammer symbol
$(c)_m =\Gamma(c+m)/\Gamma(c)$. The key point is that (\ref{discrete-Prabhakar}) 
and the entire class defined by (\ref{genfugen}) are discrete approximations of the so called
Prabhakar or generalized Mittag--Leffler density,
which is recovered 
by the well-scaled continuous-time limit \cite{Michelitsch-Polito-Riascos2021} 
\\[2mm]
$
\displaystyle \chi_{\alpha,\nu }(t)  =\lim_{h \to 0}\frac{1}{h} \varphi_{\alpha}^{(\nu)}\left(\frac{t}{h},\xi_0h^{\alpha}\right)   =\lim_{h \to 0}\frac{1}{h} \theta_{\alpha}^{(\nu)}\left(\frac{t}{h},\xi_0h^{\alpha}\right) =  \xi_0^{\nu } t^{\nu \alpha-1} E_{\alpha,\nu \alpha}^{\nu }(-\xi_0 t^{\alpha})
$\\ \\ ($\alpha \in (0,1]$,\,\, $\nu>0$, $\xi>0$,\,\, $t\in h {\mathbf N} \to {\mathbf R}^{+}$).
The multiplier $h^{-1}$ reflects the fact that this is a density (i.e.\ of physical dimension $sec^{-1}$). Note that the above expression
contains the Prabhakar function $E_{a,b}^c(z) = \sum_{m=0}^{\infty}
 \frac{(c)_m}{m!}\frac{z^m}{\Gamma(am + b)}$ 
 ($\Re\{a\} >0, \Re\{b\} > 0, c,\, z \in \mathbf{C}$) \cite{Prabhakar1971}.
 The continuous-time Prabhakar counting process
 with inter-arrival time density
 $\chi_{\alpha,\nu }(t)$ was first introduced by Cahoy and Polito as a generalization of the fractional Poisson process \cite{PolitoCahoy2013}.
 Here we call this process `{\it Generalized fractional Poisson Process (GFPP)}'. The GFPP was applied to random walks on undirected graphs by Michelitsch and Riascos \cite{TMM-APR-PhysicaA2020,MiRia2020}. We point out that the time-scale parameter $\xi(h)=\xi_0h^{\alpha}$ is re-scaled in such a way that a well-scaled continuous-time limit to the Prabhakar density exists and a new time-scale parameter $\xi_0>0$ of physical dimension $sec^{-\alpha}$ and independent of $h$ comes into play.
 By simple conditioning arguments and using the connection of generating functions with discrete-time convolutions one obtains the PDTP state probability generating function ${\bar \Phi}^{(n)}_{\alpha,\nu}(u,\xi)= \frac{(1-u{\bar \varphi}_{\alpha}^{(\nu)}(u))}{(1-u)} u^n{\bar \varphi}_{\alpha}^{(n\nu)}(u)$,
which in turn leads to the PDTP state probabilities (probabilities of having $n$ arrivals within $[0,t]$) as follows. If $t < n $, $\Phi^{(n)}_{\alpha,\nu}(t,\xi) = 0$, otherwise  
\begin{equation}
\label{state-prob}
\begin{array}{lll}
\displaystyle \Phi^{(n)}_{\alpha,\nu}(t,\xi)
  = \displaystyle
\frac{\Theta(t-n)\xi^{n\nu}}{(t-n)!}\sum_{m=0}^{\infty} \frac{(-1)^m (n\nu)_{m} \xi^m}{m!}  \frac{\Gamma[\alpha(m+n\nu)+1+t-n]}{\Gamma[\alpha(m+n\nu)+1]} & \\ \\ \displaystyle 
-\frac{\Theta(t-n-1)\xi^{(n+1)\nu}}{(t-n-1)!}\sum_{m=0}^{\infty} \frac{(-1)^m ([n+1]\nu)_{m} \xi^m}{m!}  \frac{\Gamma[\alpha[m+(n+1)\nu]+t-n]}{\Gamma(\alpha[m+(n+1)\nu]+1)} ,
 \\ \\ \mbox{for} \: 0<\xi<1, \mbox{while for } \xi >1  &
\\ \\
\displaystyle \Phi^{(n)}_{\alpha,\nu}(t,\xi)
  = \displaystyle
\frac{\Theta(t-n)(-1)^{t-n}}{(t-n)!}
\sum_{m=0}^{\infty} \frac{(-1)^m (n\nu )_{m} \xi^{-m}}{m!}
\frac{\Gamma[\alpha m]}{\Gamma[\alpha m-(t-n)]} \\ \\ \hspace{0.5cm} \displaystyle -\frac{\Theta(t-n-1)(-1)^{t-n-1}}{(t-n-1)!}
\sum_{m=0}^{\infty} \frac{(-1)^m ([n+1]\nu )_{m} \xi^{-m}}{m!}
\frac{\Gamma[\alpha m]}{\Gamma[\alpha m-(t-n-1)]} ,
\end{array} 
\end{equation}
where $n,t \in {\mathbf N_0}$ and we introduced the discrete-time version of the Heaviside step function $\Theta(r)$ with $\Theta(r)=1$ for $r\geq 1$ and $\Theta(r)=0$ otherwise. Note that $\Theta(0)=1$. We observe that
the state probabilities are non-zero only for
$0\leq n \leq t$ reflecting the general feature that there cannot be more arrivals than time increments in a discrete-time counting process. Crucial are the initial condition $\Phi^{(n)}_{\alpha,\nu}(t,\xi)\big|_{t=0} = \delta_{n0}$ (with the Kronecker symbol $\delta_{ij}$) and the normalization condition $\sum_{n=0}^t\Phi^{(n)}_{\alpha,\nu}(t,\xi)=1$.
The PDTP state probabilities (\ref{state-prob}) are an approximation of the GFPP state probabilities (the latter first obtained
by Cahoy and Polito \cite{PolitoCahoy2013} and consult also \cite{TMM-APR-PhysicaA2020,MiRia2020}). The probabilities (\ref{state-prob}) are connected to the GFPP state probabilities by the well-scaled continuous-time limit
\begin{equation}
    \label{well-sclaed_state-prob}
\label{theform}
\begin{array}{clr}
\Phi_{\alpha,\nu}^{(n)}(t,\xi_0)_{ct} & =  \lim_{h\rightarrow 0} \Phi_{\alpha,\nu}^{(n)}\left(\frac{t}{h},\xi_0h^{\alpha}\right), \hspace{0.3cm}  n\in \mathbf{N}_0, 
 \hspace{0.1cm} t \in h\mathbf{N}_0 \to {\mathbf R}^{+}  & \\ \\
  & =  (\xi_0t^{\alpha})^{n\nu} \left\{ E_{\alpha, \alpha n \nu +1}^{n\nu}(-\xi_0t^{\alpha}) -  
 (\xi_0t^{\alpha})^{\nu} E_{\alpha, \alpha (n+1) \nu +1}^{(n+1)\nu}(-\xi_0t^{\alpha})\right\}. &
\end{array}
\end{equation}
The continuous-time GFPP state probabilities are plotted in Fig. \ref{Fig1} for different parameters. One can see that higher states $n$ are `less occupied' at time $t$.
Worthy of mention is also the universal power-law long-time asymptotic feature
\begin{equation}
\label{long-time}
 \Phi_{\alpha,\nu}^{(n)}(t,\xi)  \sim \frac{\nu}{\xi} \frac{t^{-\alpha}}{\Gamma(1-\alpha)}, \hspace{1cm} \alpha \in (0,1),
\end{equation}
which is independent of state $n$, i.e.\ at large times all states have the same power-law tail behavior.
This is caused by the long-time memory and non-markovianity of the PDTP and induces a fat-tail in the inter-arrival time density  $\theta_{\alpha}^{(\nu)}(t,\xi) \sim  \frac{\alpha\nu}{\xi_0} \frac{t^{-\alpha -1}}{\Gamma(1-\alpha)}$, $\alpha \in (0,1)$, which is of the same type as for the Mittag--Leffler density.
\section{Prabhakar discrete time random walk model}
\label{Prabhakar_DTRW}
Here we consider a random walk on undirected connected graphs of $N$ nodes subordinated to a PDTP. 
The resulting walk is a discrete-time analogue to the classical Montroll--Weiss CTRW (see \cite{MontrollWeiss1965} for details). We call this walk Prabhakar Discrete-time Random Walk (Prabhakar DTRW).
The jumps of the walker from node $i \to j$ are governed by the $N\times N$ one-step transition matrix ${\mathbf H}$
with the elements 
$H_{ij} = \frac{1}{K_i}A_{ij}$, where $K_i$ and $A_{ij}$ indicate the degree and (in undirected graphs symmetric) adjacency matrix, respectively (see \cite{NohRieger2004,TMM-APR-ISTE2019,RiascosMateos2014} for details).
The $N\times N$ Prabhakar DTRW transition matrix then writes
\begin{equation}
\label{Coxseries}
{\mathbf P}^{\alpha,\nu}(t,\xi) = \sum_{n=0}^{t} {\mathbf H}^n \Phi^{(n)}_{\alpha,\nu}(t,\xi) ,\hspace{1cm} 
P_{ij}(t)\big|_{t=0}= \delta_{ij} ,\hspace{1cm} t \in \mathbf{N}_0,
\end{equation}
where this series stops at $n=t$ due to the
above discussed properties of the state probabilities. Thus, (\ref{Coxseries}) is a matrix polynomial in ${\mathbf H}$ of degree $t$. The matrix element $P_{ij}^{\alpha,\nu}(t,\xi)$ indicates the probability that the walker is present on node $j$ at time $t$ when having started the walk at $t=0$
on node $i$. The following generalized fractional difference equation holds:
\begin{equation}
\label{straight}
 {\hat {\cal D}}_{\alpha}^{\nu} \cdot {\mathbf P}^{\alpha,\nu}(t,\xi) - {\cal M}_{\alpha,\nu}(t,\xi) 
{\mathbf 1}
= {\mathbf H} \cdot {\mathbf P}^{\alpha,\nu}(t-1,\xi)   ,\hspace{1cm} P^{\alpha,\nu}_{ij}(0) =\delta_{ij}.
\end{equation}
In the above equation, ${\cal M}_{\alpha,\nu}(t,\xi)$ is a memory kernel
and ${\hat {\cal D}}_{\alpha}^{\nu}$ denotes a Prabhakar generalized fractional difference operator, see \cite{Michelitsch-Polito-Riascos2021} for their explicit forms and further details. The well-scaled continuous-time limit of the transition matrix (\ref{Coxseries}), defined by ${\mathbf P}^{\alpha,\nu}(t,\xi_0)_{ct}=
\lim_{h\to 0}{\mathbf P}^{\alpha,\nu}(\frac{t}{h},\xi_0h^{\alpha})$, yields the transition matrix of the walk subordinated to the GFPP which is a Montroll--Weiss CTRW governed by a generalized fractional Kolmogorov-Feller equation involving Prabhakar kernels (see the review article \cite{GiustiPolitoMainardi-etal2020} as well as \cite{Michelitsch-Polito-Riascos2021} and consult \cite{TMM-APR-PhysicaA2020,MiRia2020} for the CTRW model).
\section{Conclusions}
Discrete-time counting processes such as
the PDTP and the related DTRW models have large potential of new interdisciplinary
applications to topical problems, including generalizations of existing models such as in biased random walks \cite{RiascosMichelitschPisarroM2020}, the time evolution of pandemic spread \cite{BesRiaMiCo2021}, communications in complex
networks, dynamics in public transport networks, anomalous relaxation and transport, population dynamics, collapse of financial markets, just to name a few, and generally to the modelling of `complex systems'.
%
%

%

\end{document}